\documentclass{article}

\usepackage[utf8]{inputenc}
\usepackage[english]{babel}

\usepackage{amsthm}
\usepackage{graphicx}
\usepackage{amssymb}
\usepackage{amsmath}
\usepackage{tikz-cd}
\usepackage{enumerate}
\usepackage{verbatim}
\usepackage{mathtools}
\usepackage{geometry}

\DeclareRobustCommand{\qbinom}{\genfrac{\lbrack}{\rbrack}{0pt}{}}
\def\multiset#1#2{\ensuremath{\left(\kern-.3em\left(\genfrac{}{}{0pt}{}{#1}{#2}\right)\kern-.3em\right)}}

\newtheorem*{theorem*}{Theorem}

\begin{document}
\title{The ``Shape" of $q$-Binomial Coefficients: \\
\vspace{.3cm}
\small{A friendly Q \& A.} }
\author{by Nate Harman}

%\affil{ Massachusetts Institute of Technology}

%\email{nharman@math.mit.edu}

%\subjclass{Primary 05A16}

\newgeometry{left=3cm,bottom=0.1cm}

\date{}
\maketitle

\begin{center}
\includegraphics[width=50mm]{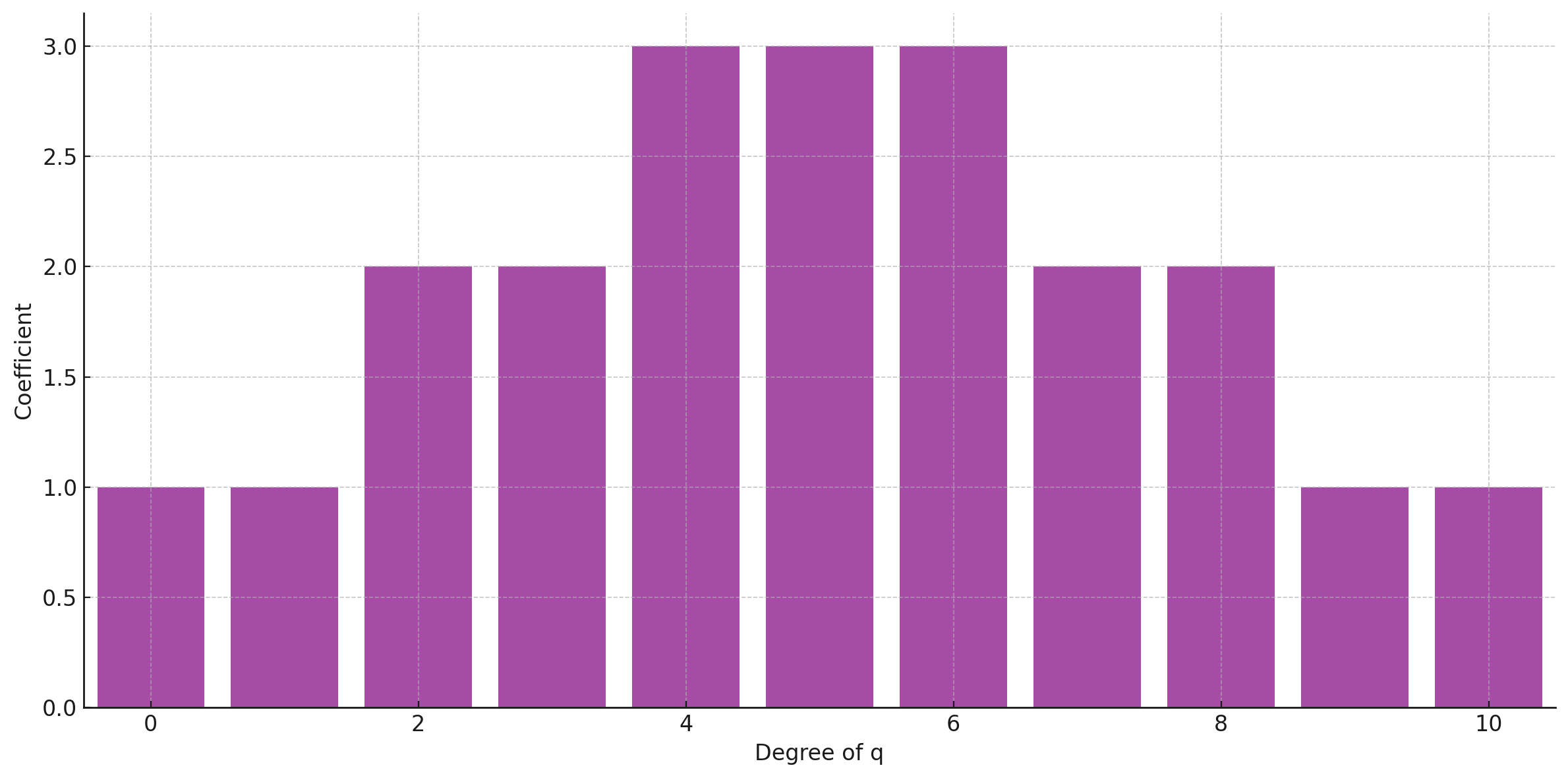}
\includegraphics[width=50mm]{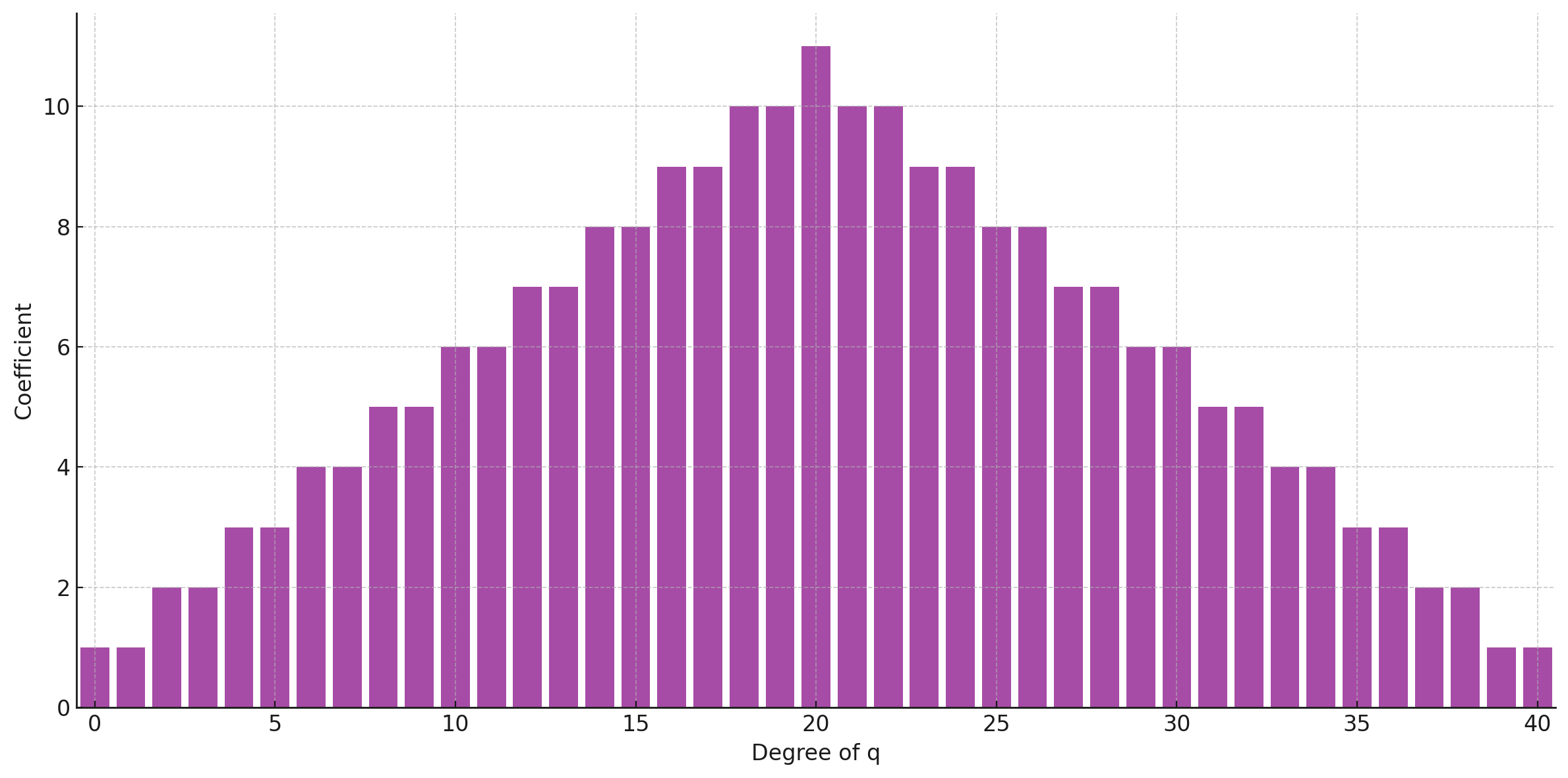}
\includegraphics[width=50mm]{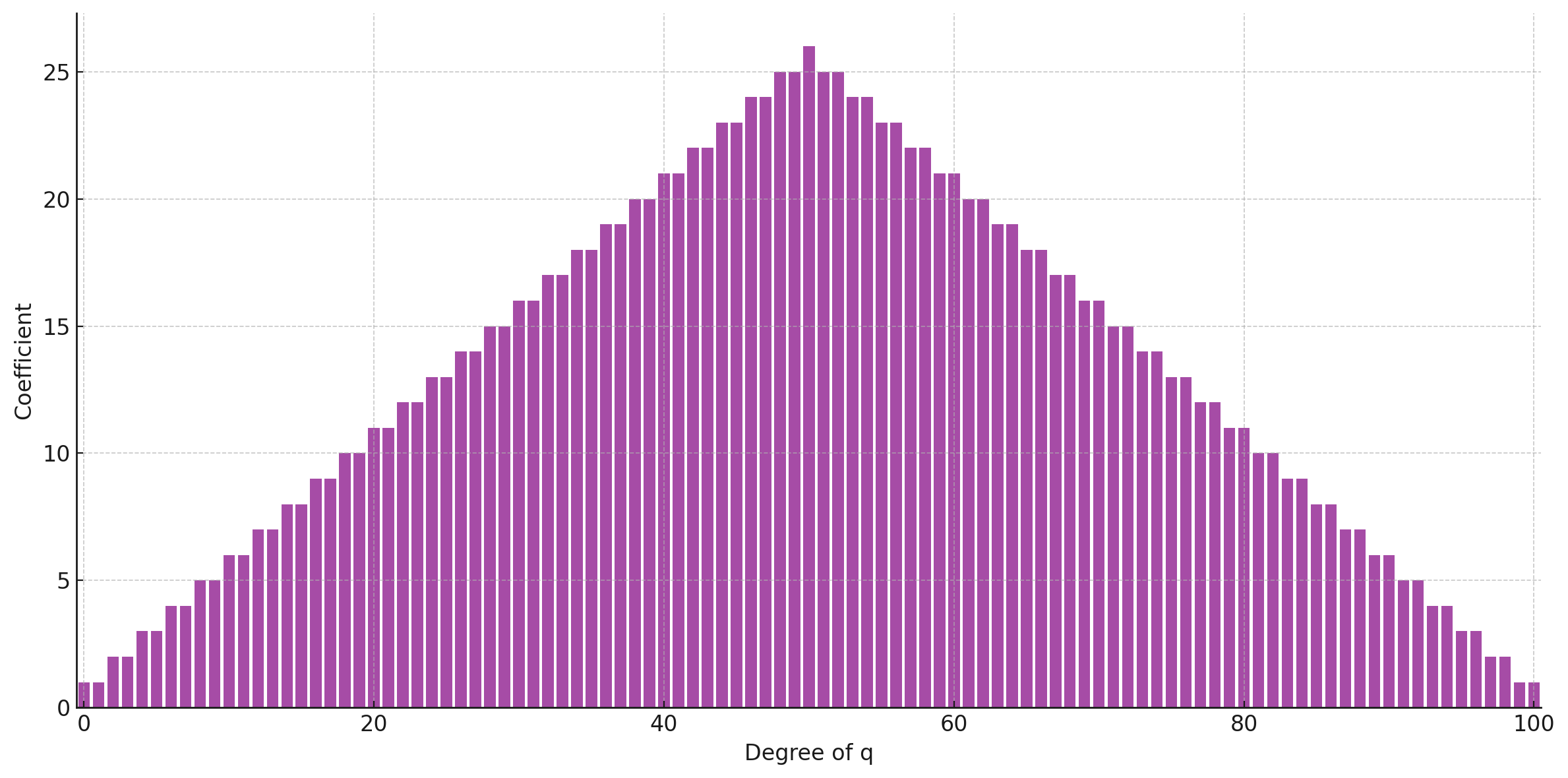}
\\
The coefficients of $\qbinom{n+2}{2}_q$ for $n = 5$, $n=20$, and $n=50$.
\\
\includegraphics[width=50mm]{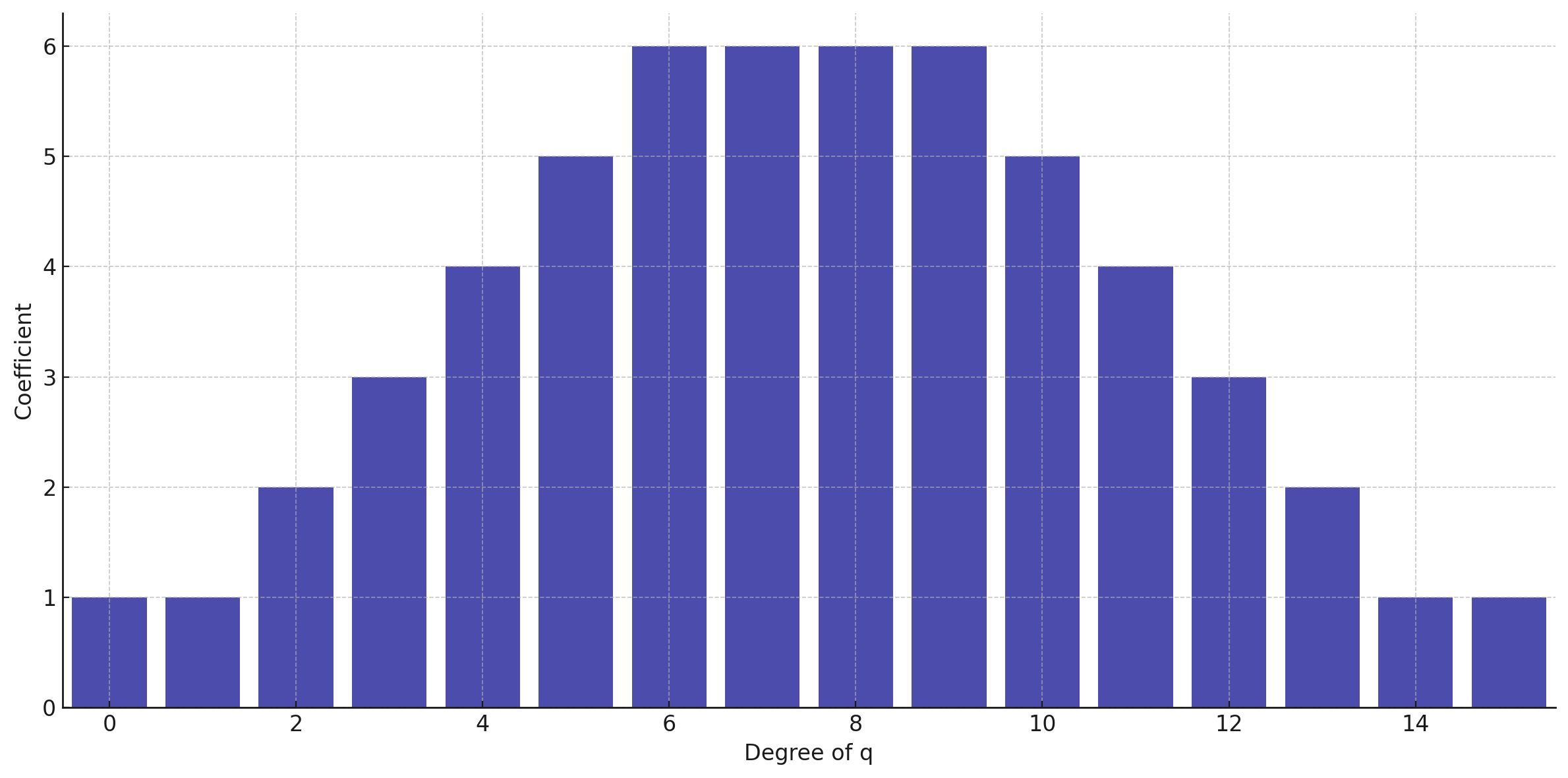}
\includegraphics[width=50mm]{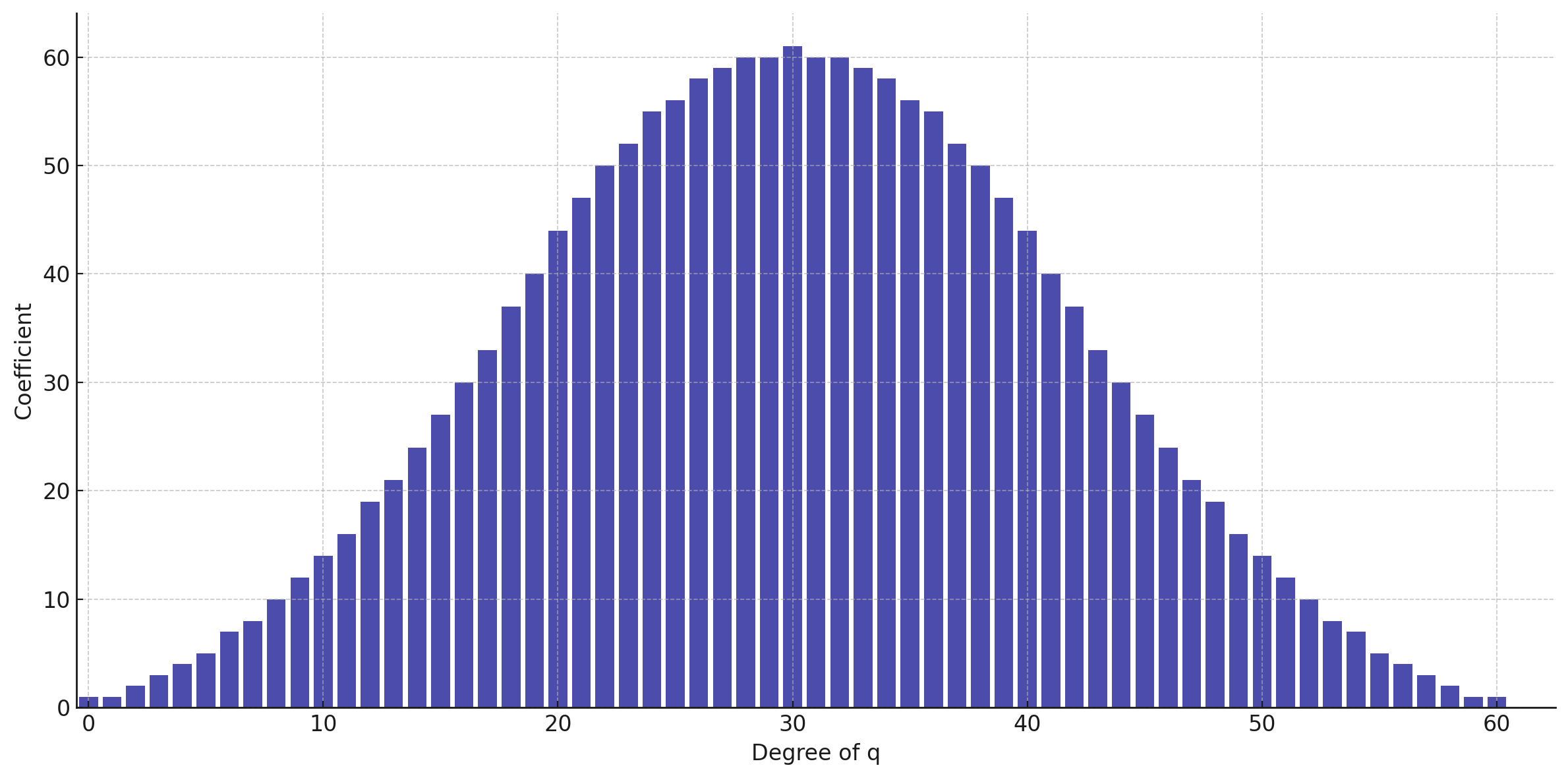}
\includegraphics[width=50mm]{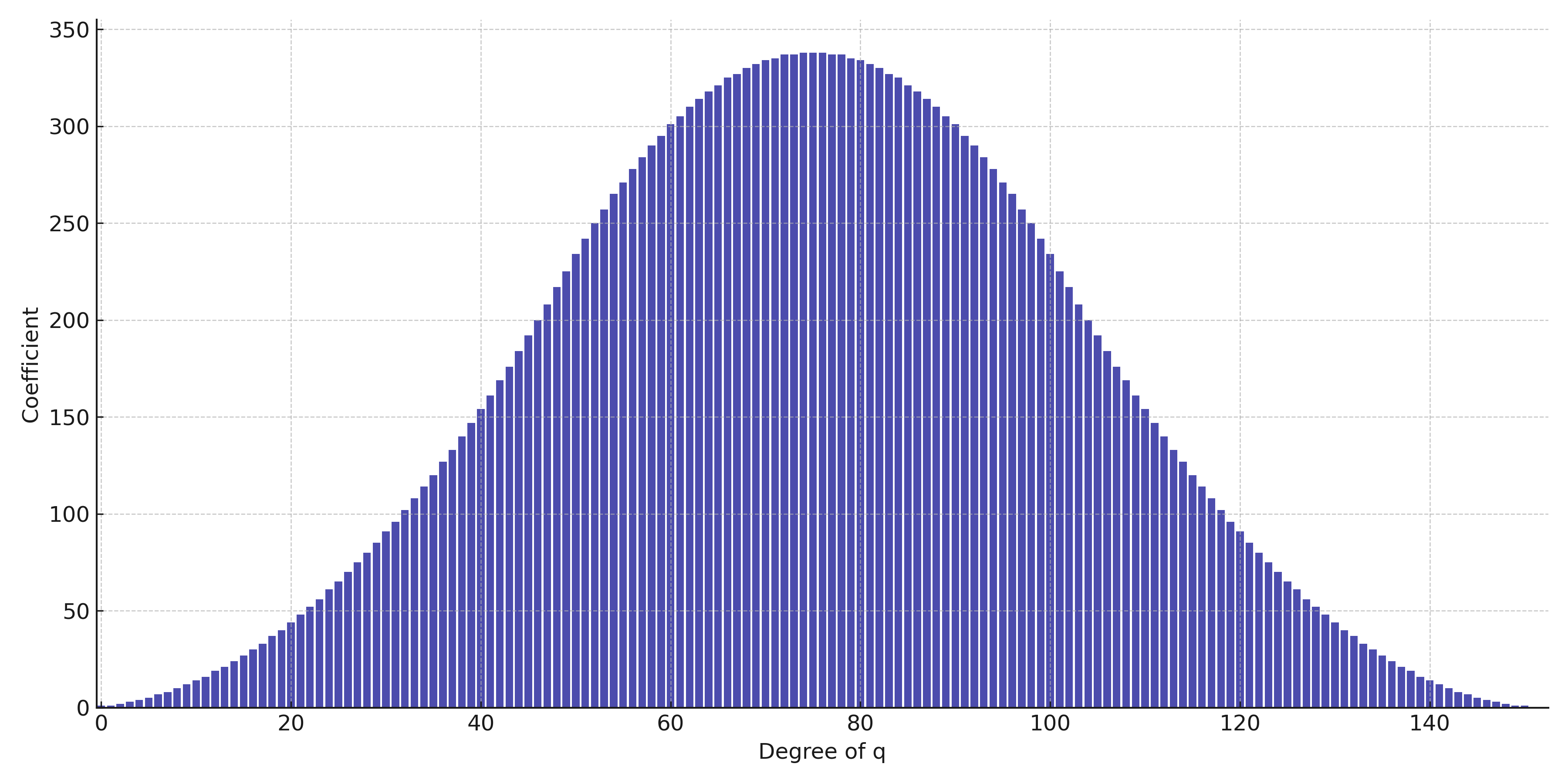}
\\
The coefficients of $\qbinom{n+3}{3}_q$ for $n = 5$, $n=20$, and $n=50$.
\\
\includegraphics[width=50mm]{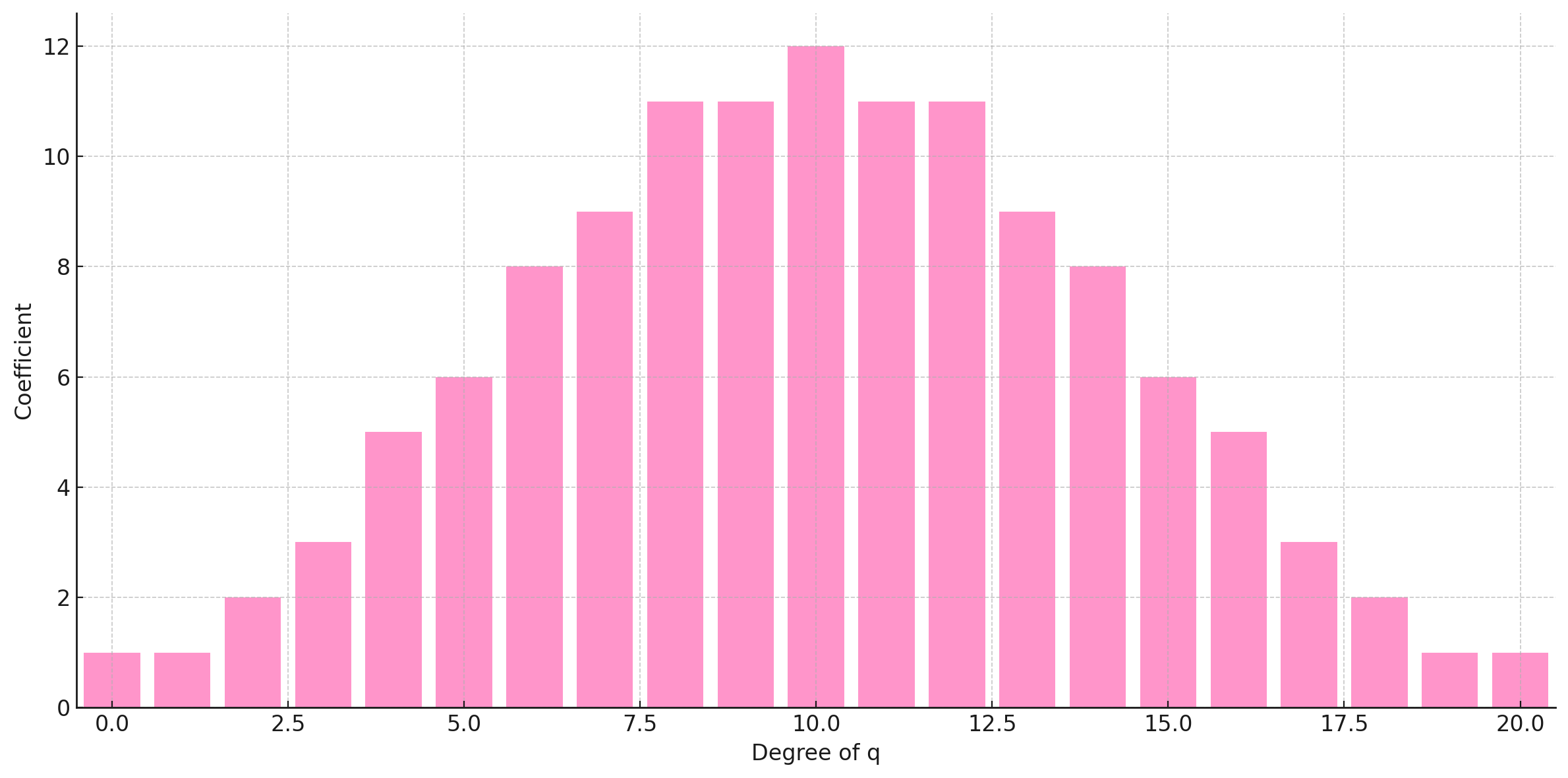}
\includegraphics[width=50mm]{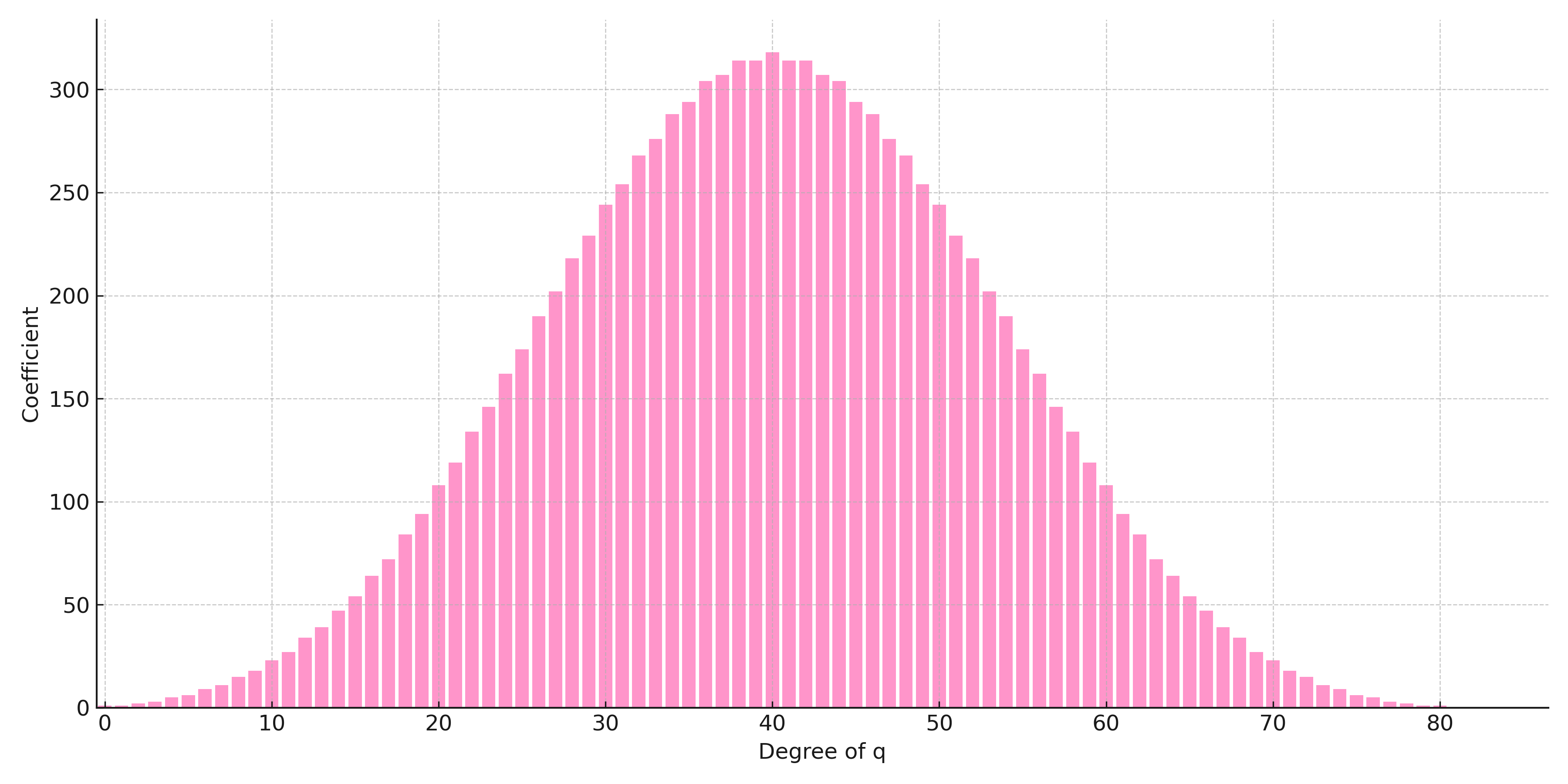}
\includegraphics[width=50mm]{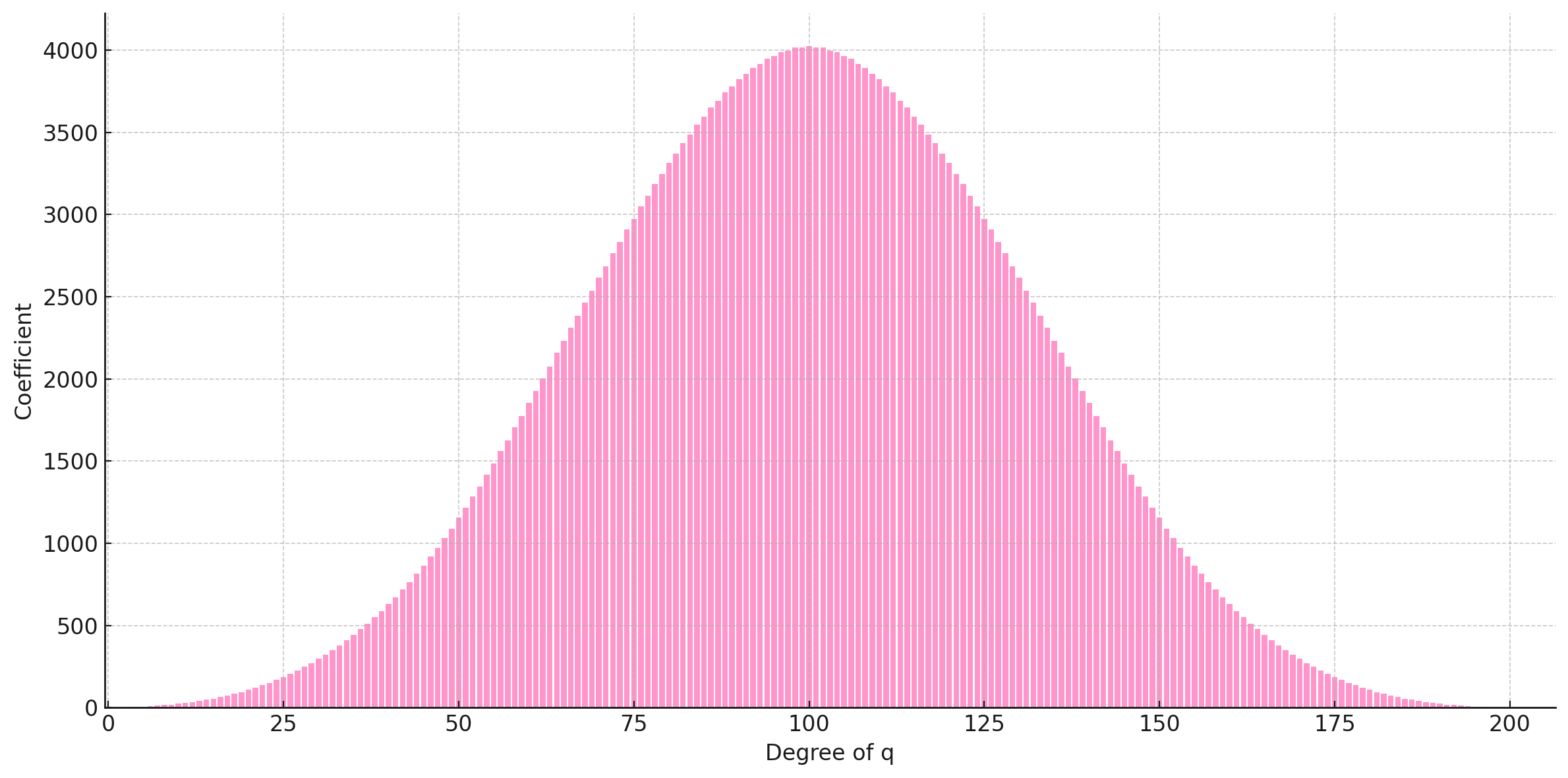}
\\
The coefficients of $\qbinom{n+4}{4}_q$ for $n = 5$, $n=20$, and $n=50$.

\end{center}

\vspace{1cm}

\begin{abstract}
Look at those shapes! What's up with that?

\end{abstract}

\restoregeometry

\noindent \textbf{Q. Why are you writing an expository Q \& A?}

\medskip

\noindent A.  All of the actual facts and theorems in this note have already occurred in the literature.  What was missing (in my opinion) was the pictures:  both the literal pictures I included to demonstrate these ideas, as well as the ``big picture" story of how some of these ideas and results fit together.  

 I have known the story presented here, and had these pictures in my head for some time now.  The real impetus for me writing this was just me playing with ChatGPT and getting it to generate some of these nice diagrams.   Once I had the diagrams I felt like I might as well share them.

\medskip

\noindent The Q \& A style is just for fun.

\bigskip

\noindent \textbf{Q. What are $q$-binomial coefficients and why do we care?}

\medskip

\noindent A. First we define the $q$-integers:
$$[n]_q = \frac{1-q^n}{1-q} =  1 + q +q^2 + \dots + q^{n-1}$$
then we define the $q$-factorials:
$$[n]!_q = [n]_q [n-1]_q [n-2]_q \dots [2]_q [1]_q$$
and finally we define the $q$-binomial coefficients: 
$$\qbinom{n}{k}_q = \frac{[n]!_q}{[n-k]!_q [k]!_q} =  \frac{(1 - q^{n+k})(1 - q^{n+k-1})\cdots(1 - q^{n+1})}{(1-q^k)(1-q^{k-1}) \cdots (1-q^2)(1-q)}. $$

\noindent  It is not immediately obvious from the formula but the $q$-binomial coefficient $\qbinom{n}{k}_q$ is in fact a polynomial in $q$. If we specialize to $q=1$ we recover the usual binomial coefficient $\binom{n}{k}$.  As for why we care, here are a few possible answers:

\begin{itemize} 

\item Combinatorially $\qbinom{n+k}{k}_q$ can be thought of as a generating function $\qbinom{n+k}{k}_q = \sum_{i=0}^{nk} P(i,n,k)q^i$ where $P(i,n,k)$ counts the number of partitions of size $i$ with at most $n$ parts, and all parts having size at most $k$.  

\item  The $q$-binomial coefficient $\qbinom{n}{k}_q$ is the Hilbert-Poincar\'e series for the complex Grassmannian $$\qbinom{n}{k}_q = \sum_{i=0}^{nk} dim(H^{2i}(Gr(n,k, \mathbb{C}),\mathbb{Q})q^i.$$ 

\item If $q$ is a prime power, then $\qbinom{n}{k}_q = |Gr(n,k, \mathbb{F}_q)|$, the number of $k$-dimensional subspaces of $\mathbb{F}_q^n$.  (This is connected to the previous point via the Grothendieck-Lefschetz trace formula.)

\item $\qbinom{n+k}{k}_q$ is the character of the $GL_2(\mathbb{C})$-representation $Sym^k(Sym^{n}(\mathbb{C}^2))$ evaluated at the matrix $\begin{bmatrix}
1 & 0 \\
0 & q 
\end{bmatrix}
\in GL_2(\mathbb{C})$.  Such representations are known as \emph{plethysms}, they are a ubiquitous and mysterious class of representations.

\end{itemize}

\bigskip

\pagebreak

\noindent \textbf{Q. What are those pictures on the first page?}

\medskip

\noindent A. Given a polynomial in a variable $q$  we can make a bar graph where the height of $k$th bar corresponds to the coefficient of $q^k$.  For example if $f(q) = 1 + 3q +9q^2 + 4q^4 +2q^5 + 9q^6$ then our bar graph would be:  

\begin{center}
\includegraphics[width=50mm]{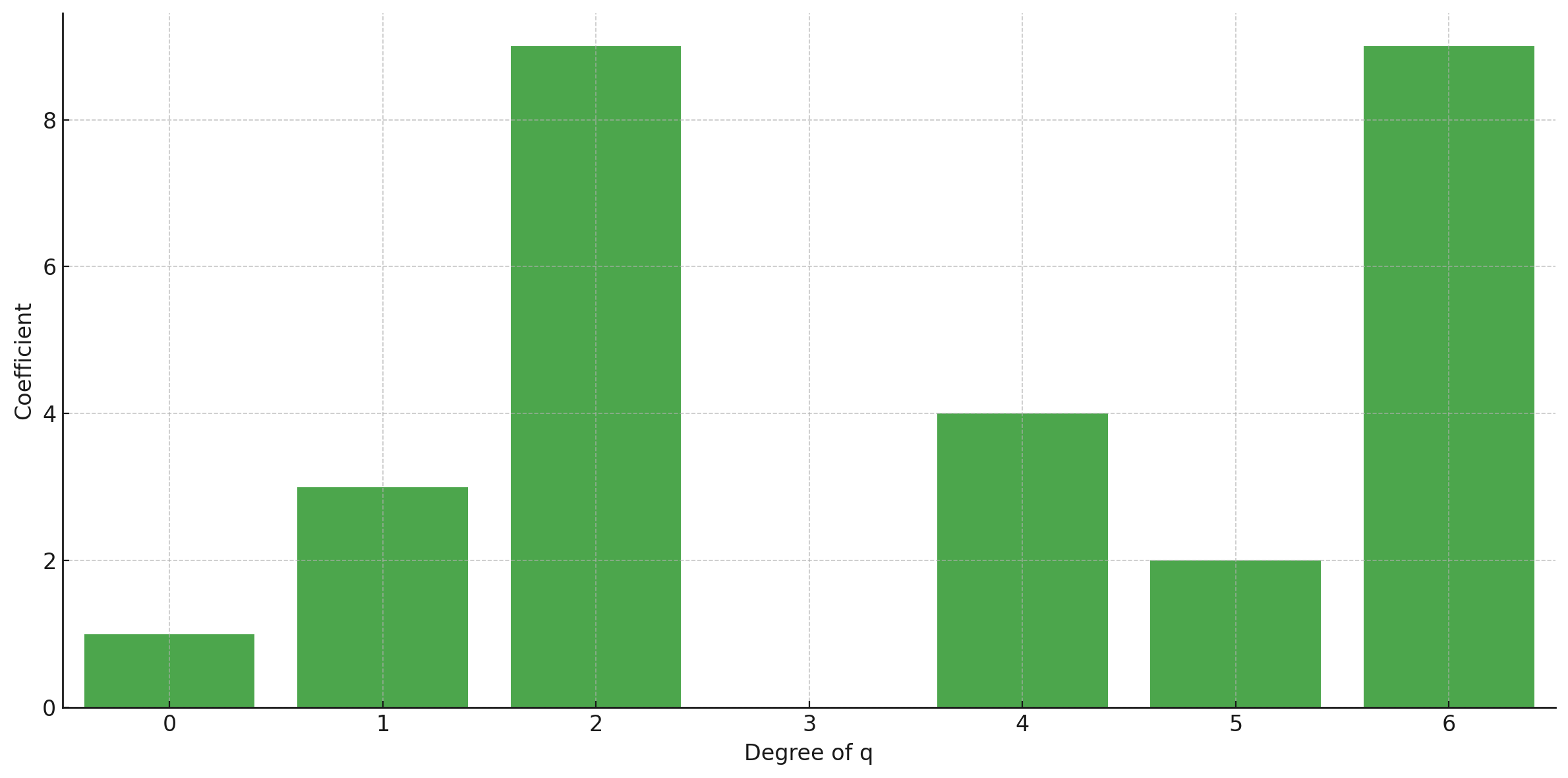}
\end{center}
 We normalize things so that the width and height of our bar graphs are constant, allowing us to compare the overall ``shapes" of polynomials of different degrees.   We are  interested in doing this for the binomial coefficients $\qbinom{n+k}{k}_q$ with $k$ fixed and $n$ growing.  Here are some pictures again for $k=3$:

\medskip

\hspace{-.5cm}
\includegraphics[width=50mm]{chart_5+3,3.png}
\includegraphics[width=50mm]{chart_20+3,3.png}
\includegraphics[width=50mm]{chart_50+3,3.png}
We can see in these pictures the well-known fact the coefficients of these polynomials are symmetric and unimodal. Less well-known though, and the focus of this note, is observation is that these bar graphs for $\qbinom{n+k}{k}_q$ are converging to a limit shape as $n \to \infty$.

\bigskip

\noindent \textbf{Q. Is there a more precise, mathematical description of what you are doing?}

\medskip

\noindent A. Given a polynomial $p(q) = a_0 + a_1q + a_2q^2 + \dots + a_n q^n$ with non-negative real coefficients,  we can construct a probability measure on $[0,1]$ by putting point masses at $\frac{k}{n}$ with measure $\frac{a_k}{a_0 + a_1 + \dots + a_n}$.  When we say a sequence of bar graphs approaches a limit shape, we really mean that this sequence of measures weakly converges to a continuous limiting distribution.

\bigskip

\noindent \textbf{Q. Are these ``limit shapes" just Gaussians?} 

\medskip 

\noindent A. A reasonable guess, but nope!  If we look at the $q$-binomial coefficient $\qbinom{n+k}{k}_q$ and let both $n$ and $k$ tend to $\infty$, then indeed it will look like a Gaussian (you should normalize slightly differently too).  However if we fix $k$ and just let $n$ grow we get a different limit shape for each $k$.  

\bigskip

\noindent \textbf{Q. Okay, so what are these limit shapes?} 

\medskip 

\noindent A.  I'm so glad you asked. The limit shape for $\qbinom{n+k}{k}_q$ as $n \to \infty$ is described by a continuous, piecewise polynomial function on $[0,1]$ which is given by a polynomial of degree $k-1$ on each of the intervals $[\frac{i}{k}, \frac{i+1}{k}]$.  For example when $k = 3$ this function is:

\[
L_3(x) = \begin{cases} \frac{27}{2}x^2 & \text{if } x \in [0, \frac{1}{3}], \\  -27x^2 +27 x - \frac{9}{2} & \text{if } x \in [\frac{1}{3}, \frac{2}{3}],\\ \frac{27}{2}(1-x)^2 & \text{if } x \in [\frac{2}{3}, 1]. \end{cases}
\]
More specifically, the limiting distribution $L_k(x)$ is the uniform distribution on an interval convolved with itself $k$-times, scaled and renormalized to be a probability distribution on $[0,1]$.

 More geometrically, $L_k$ is the $(k-1)$-dimensional volume of the intersection of the $k$-dimensional hypercube $[0,1]^k$ with the hyperplane $x_1 + x_2 + \dots + x_k = \alpha$ (again suitably renormalized).

\begin{center}
\includegraphics[width=130mm]{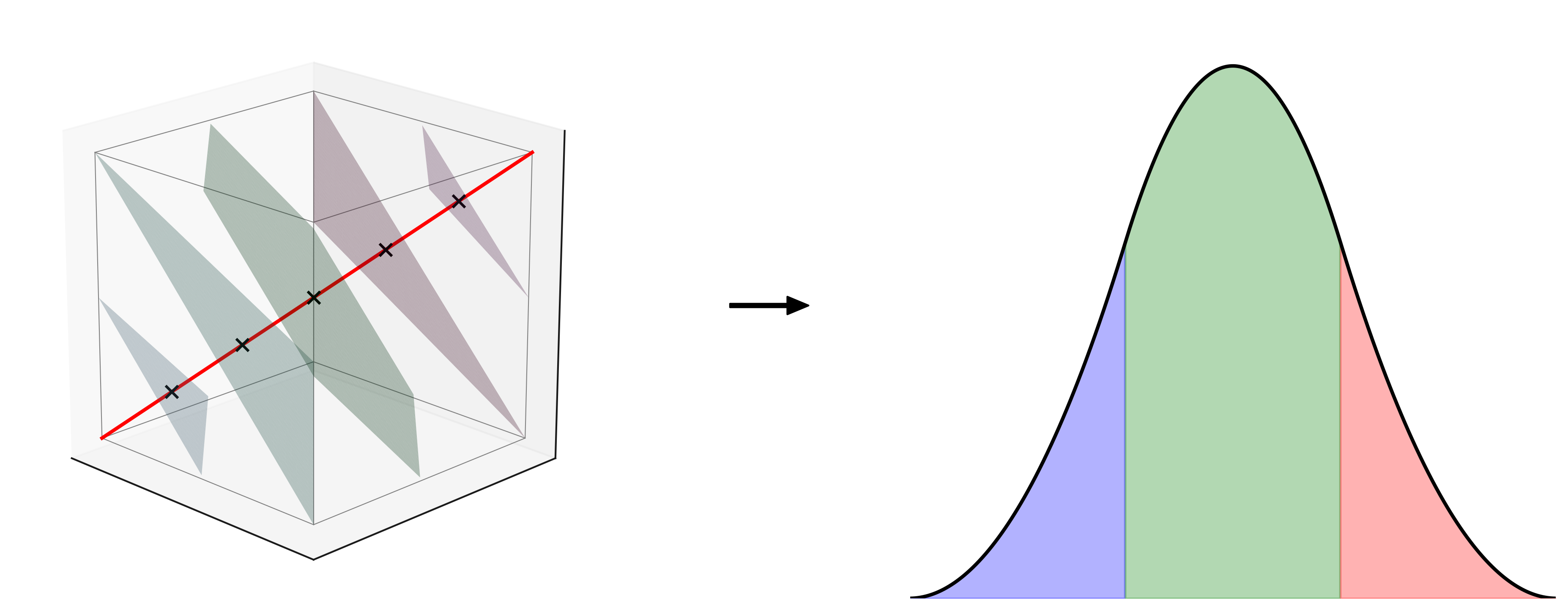}
\\
\end{center}
The places where the function switches behavior are when the hyperplane passes through a vertex and the combinatorial structure of the cross-section changes. 

\bigskip

\noindent \textbf{Q. Why a hypercube?} 

\medskip 

\noindent A.  In some sense shouldn't look at a hypercube at all. From the representation theory perspective on $q$-binomial coefficients a more natural thing to do is to look at cross sections of the simplex defined by $0 \le  x_1 \le x_2 \le \dots \le x_k \le 1$.  This is a scaled version of the Gelfand-Tsetlin polytope for $Sym^k(\mathbb{C}^{n+1})$.

  We'll note though that this simplex takes up $\frac{1}{k!}$ of the volume of the cube, and moreover its intersection with any hyperplane $x_1 +x_2 + \dots +x_k = \alpha$ is $\frac{1}{k!}$ of the $(k-1)$-dimensional volume of intersection with the cube.  So in particular, after renormalizing this gives the same measure as the full cube.

\bigskip

\noindent \textbf{Q. If the limit shapes are piecewise polynomial, can we say something similar about the finite $q$-binomials?} 

\medskip 

\noindent A. The coefficients of a $q$-binomial form an \emph{almost-piecewise} quasipolynomial function.  Don't know what almost-piecewise means?  That's because I just made it up.  Here is a color coded bar graph for $\qbinom{50+4}{4}_q$ to demonstrate what I mean:

\begin{center}
\includegraphics[width=130mm]{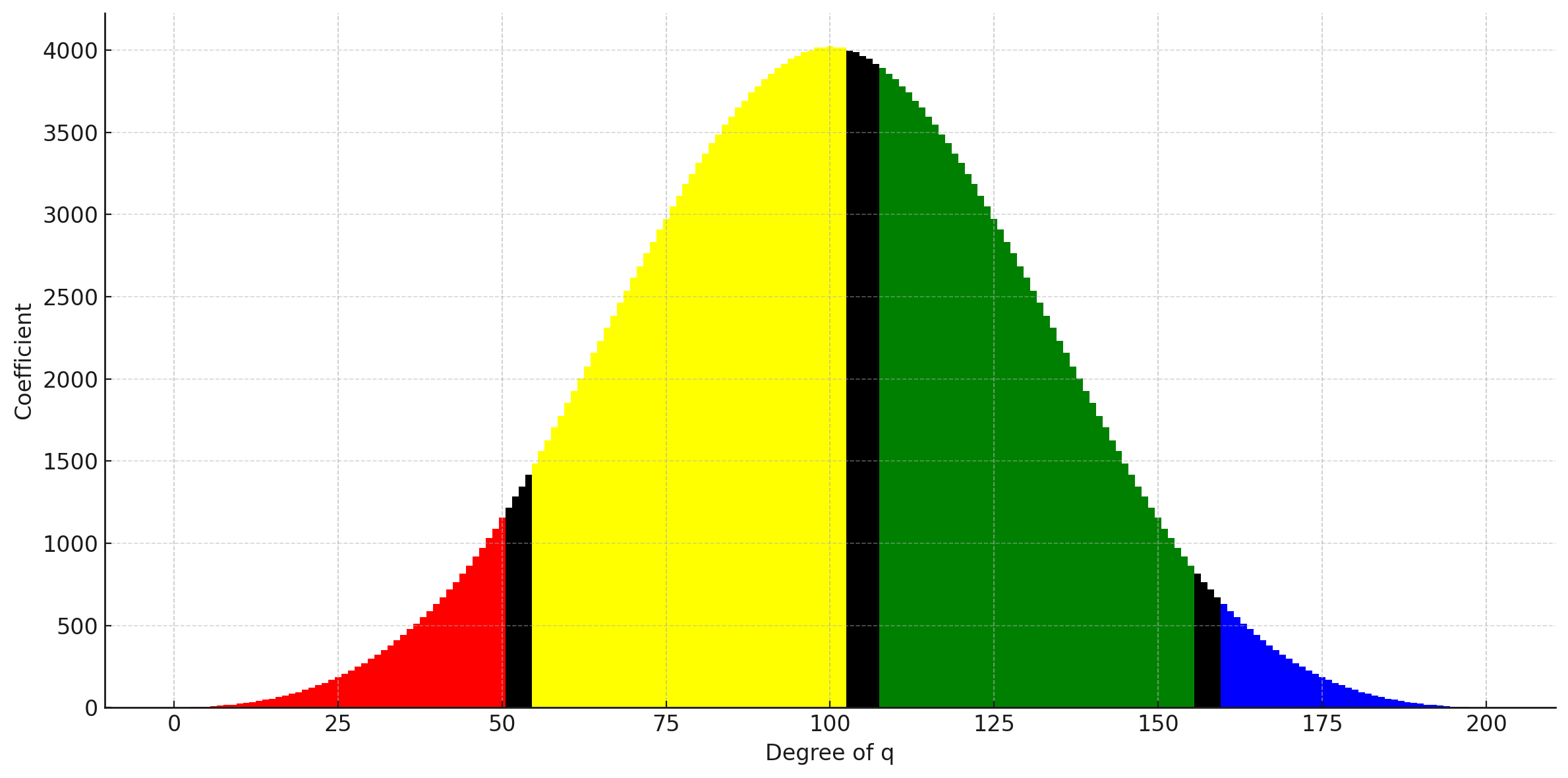}
\\
\end{center}

\begin{itemize}

\item On each of the 4 colored regions the coefficients of $\qbinom{n+4}{4}_q$ are described by a degree 3 quasipolynomial,  which is polynomial on each congruence class modulo $12$.

\item For $\qbinom{n+k}{k}_q$ we would have $k$ quasipolynomial regions, each covering about $\frac{1}{k}$ of the domain.  The quasipolynomials have degree $k-1$ and the quasiperiod is at most the l.c.m. of the numbers from $1$ to $k$. 

\item In between the colored regions we have some transition zones (colored black).  The number of bars in these zones for $\qbinom{n+k}{k}_q$ depends on $k$ but not on $n$, so they vanish in the limit.

\end{itemize}

\noindent Why is this true?  Let's look at the formula for $\qbinom{n+4}{4}_q$:

$$ \qbinom{n+4}{4}_q = \frac{(1 - q^{n+4})(1 - q^{n+3})(1 - q^{n+2})(1 - q^{n+1})}{(1-q^4)(1-q^3)(1-q^2)(1-q)}  $$ 

\noindent If we expand out the numerator this becomes:

$$ \frac{1 -  q^{n+1} - q^{n+2} - q^{n+3} - q^{n+4} + q^{2n+3} + q^{2n+4}  + 2q^{2n+5} + q^{2n+6} + q^{2n+7} - q^{3n+6} - q^{3n+7} - q^{3n+8} - q^{3n+9} + q^{4n+10}}{(1-q^4)(1-q^3)(1-q^2)(1-q)}$$ 

\noindent Notice that if we are interested in the coefficient of $q^m$ with $m \le n$ then the only term in the numerator we need to worry about is the $1$.  So for those initial values the coefficient of $q^m$ is the same as for:
$$\frac{1}{(1-q^4)(1-q^3)(1-q^2)(1-q)}$$
and those coefficients for $q^k$ can easily be seen to be quasipolynomial in $m$.  This corresponds to the red region in the picture above. Let's denote this quasipolynomial for the coefficient of $q^m$ by $f_{red}(m)$.  We will look at the explicit formulas for $f_{red}(m)$ on page $7$, but for now we don't need them.

Next up, let's think about the coefficient of $q^m$ where $n < m < 2n+3$ -- the yellow region in the picture.  Here now we do have to consider some powers of $q$ in the numerator, but not any of degree at least $2n+3$.  So for these values the coefficients agree with:
$$ \frac{1 -  q^{n+1} - q^{n+2} - q^{n+3} - q^{n+4}}{(1-q^4)(1-q^3)(1-q^2)(1-q)}$$ 
Again the theory of generating functions tells us that the coefficient of $q^m$ in this is quasipolynomial for $m \ge n+4$. Let's call this quasipolynomial $f_{yel}(x)$.  Looking at the above expression we can see directly that:
$$f_{yel}(m) = f_{red}(m) - f_{red}(m-n-1) - f_{red}(m-n-2) - f_{red}(m-n-3)  - f_{red}(m-n-4)$$ 
for $m \ge n+4$.  So if we have quasipolynomial expression for $f_{red}(m)$, we can get one for $f_{yel}(m)$ as well using this recursion. 

The black transition zone between the red and yellow regions occurs where the $q^{n+j}$ terms first appear.  For those we only get some of the terms in the above recursion, for example if we want the coefficient of $q^{n+2}$ we only need to subtract off $2$ terms instead of $4$:
$$ [q^{n+2}] = f_{red}(n+2) - f_{red}(1) - f_{red}(0) $$

We could keep going from here.  If we want a quasipolynomial expression on the green region we need to also include the $q^{2n + j}$ terms in the numerator.  Similarly for the blue region we need to include the $q^{3n + j}$ terms as well. In both cases we can express the quasipolynomial functions $f_{grn}(m)$ and $f_{blu}(m)$ in terms of the original quasipolynomial function $f_{red}(m)$.

\medskip
\noindent \textbf{Q. These ``transition zones" are scary. Do we really need them?} 

\smallskip 

\noindent A. Probably not.  The analysis of the numerator above does seem to give us tight information about when each quasipolynomial formula stops holding -- the right endpoints of the colored regions.  However we were not particularly careful about where they start holding -- the left endpoints of the colored regions.   Based on symmetry, we expect for the picture actually looks like: 
\vspace{-.2cm}
\begin{center}
\includegraphics[width=127mm]{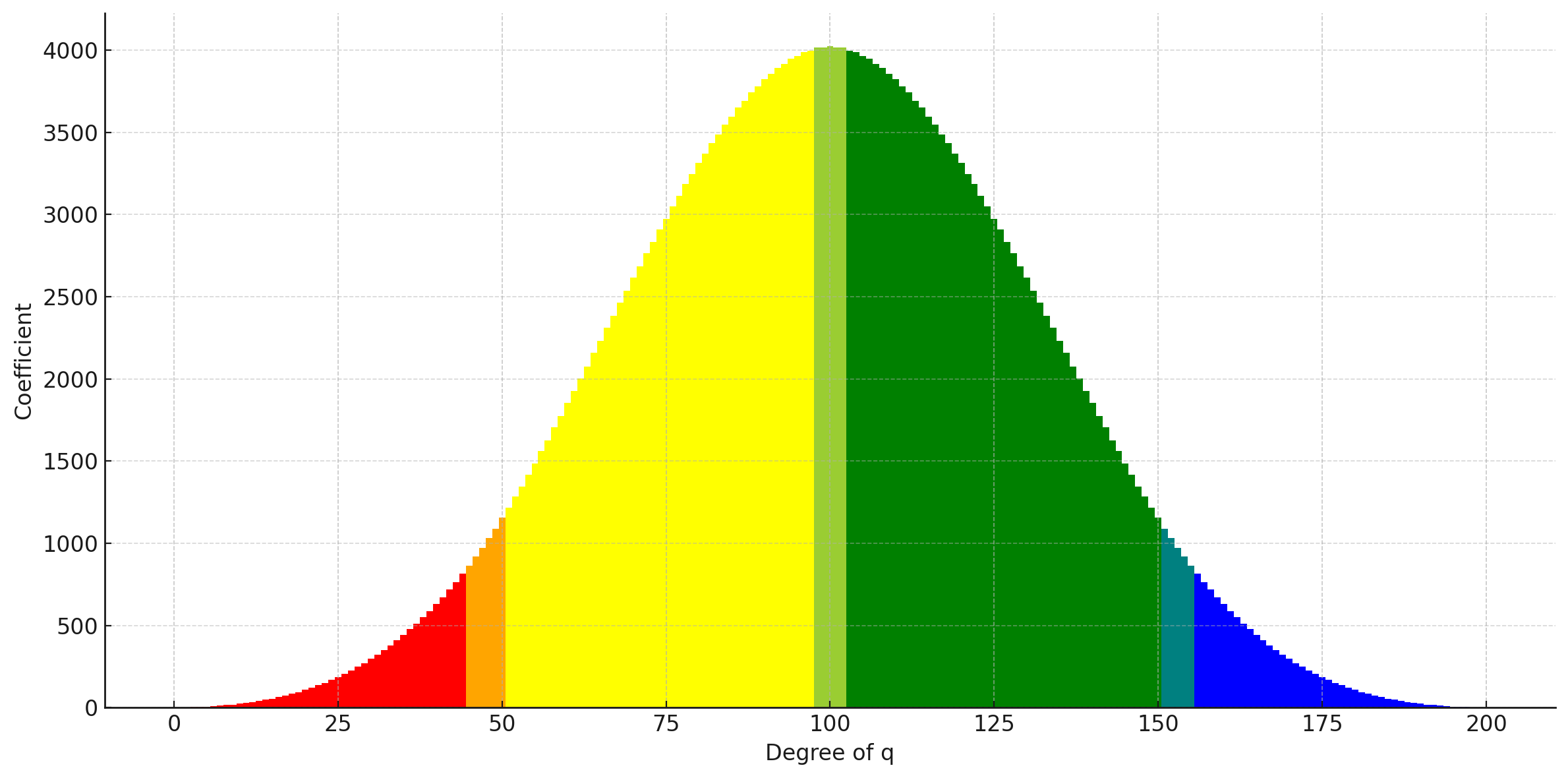}
\\
\end{center}
where the quasipolynomial regions actually overlap with one another, and agree on the overlap.

\bigskip
\noindent \textbf{Q. Is quasipolynomiality what makes these bar graphs look so nice and smooth?} 

\medskip 

\noindent A. Not quite on its own.  Let's look at the bar graph for the quasipolynomial function

$$f(x) = \begin{cases} 10x & \text{if } x \text{ is even} ,\\ \frac{1}{2}(x^2-x) & \text{if } x  \text{ is odd} \end{cases}$$

\begin{center}
\includegraphics[width=70mm]{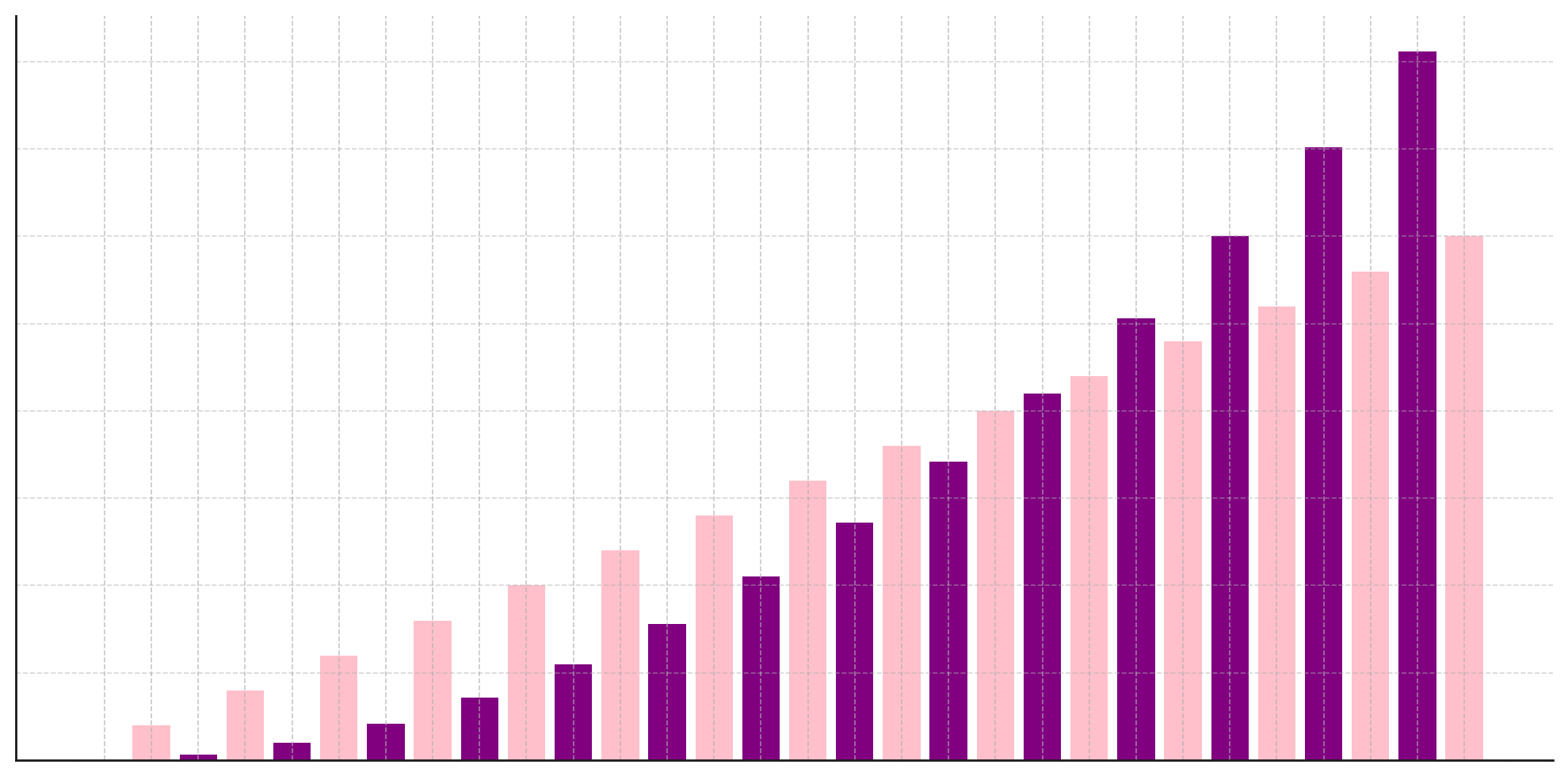}
\\
\end{center}
The behavior of the function on the even and odd values is qualitatively different, and you can see in the picture that they do not mesh together to form one nice smooth curve.  Moreover we'll note that it possible for a sequence of quasipolynomials to look like this locally but still converge in measure to a continuous distribution.  So quasipolynomiality alone can't explain why the graphs for $\qbinom{n+k}{k}_q$ look as nice as they do for finite values of $n$.  

To see why our pictures look so nice, let's look more carefully at the quasipolynomial $f_{red}(m)$ which describes the coefficient of $q^m$ in $\qbinom{n+4}{4}_q$ for $m < n$:  

\[
f_{red}(m) =
\begin{cases}
\frac{1}{144}m^3 + \frac{5}{48}m^2 + \frac{1}{2}m + 1, & m \equiv 0 \pmod{12} \\
\frac{1}{144}m^3 + \frac{5}{48}m^2 + \frac{7}{16}m + \frac{65}{144}, & m \equiv 1 \pmod{12} \\
\frac{1}{144}m^3 + \frac{5}{48}m^2 + \frac{1}{2}m + \frac{19}{36}, & m \equiv 2 \pmod{12} \\
\frac{1}{144}m^3 + \frac{5}{48}m^2 + \frac{7}{16}m + \frac{9}{16}, & m \equiv 3 \pmod{12} \\
\frac{1}{144}m^3 + \frac{5}{48}m^2 + \frac{1}{2}m + \frac{8}{9}, & m \equiv 4 \pmod{12} \\
\frac{1}{144}m^3 + \frac{5}{48}m^2 + \frac{7}{16}m + \frac{49}{144}, & m \equiv 5 \pmod{12} \\
\frac{1}{144}m^3 + \frac{5}{48}m^2 + \frac{1}{2}m + \frac{3}{4}, & m \equiv 6 \pmod{12} \\
\frac{1}{144}m^3 + \frac{5}{48}m^2 + \frac{7}{16}m + \frac{65}{144}, & m \equiv 7 \pmod{12} \\
\frac{1}{144}m^3 + \frac{5}{48}m^2 + \frac{1}{2}m + \frac{7}{9}, & m \equiv 8 \pmod{12} \\
\frac{1}{144}m^3 + \frac{5}{48}m^2 + \frac{7}{16}m + \frac{9}{16}, & m \equiv 9 \pmod{12} \\
\frac{1}{144}m^3 + \frac{5}{48}m^2 + \frac{1}{2}m + \frac{23}{36}, & m \equiv 10 \pmod{12} \\
\frac{1}{144}m^3 + \frac{5}{48}m^2 + \frac{7}{16}m + \frac{49}{144}, & m \equiv 11 \pmod{12}
\end{cases}
\]
Notice that all of the degree $3$ and degree $2$ terms are the same, and only the lower order terms have periodic coefficients. 

\begin{itemize}

\item Because we can write them in terms of values of $f_{red}$, the other quasipolynomials $f_{yel}$, $f_{grn}$, and $f_{blu}$ describing the coefficients of $\qbinom{n+4}{4}_q$  on different intervals also have this property.

\item The consecutive differences $\Delta f_{red}(m) = f_{red}(m) - f_{red}(m-1)$ are quasipolynomial in $m$ of degree $2$, with the leading quadratic terms being the same on all congruence classes.

\item More generally, if we look at the coefficients $\qbinom{n+k}{k}_q$  then the quasipolynomials have degree $k-1$ but the smallest periodic coefficient has degree $\lfloor \frac{k}{2} \rfloor -1$.    This can be easily seen for the initial terms (the analog of $f_{red}(m)$), which are coming from
$$\frac{1}{(1-q)(1-q^2)\cdots(1-q^k)}$$

\end{itemize}

\noindent This means that on a large scale, these functions really ``look like" a single polynomial and the deviation from actually fitting a single polynomial is happening on a much smaller scale.  Moreover, if $k \ge 4$ then the consecutive differences are well-approximated by the derivative of the polynomial main term.  We will take this as a moral explanation for why these bar graphs for $\qbinom{n+k}{k}_q$ look so smooth, even for relatively small values of $n$.

\bigskip
\noindent \textbf{Q. Where is this stuff actually proven?} 

\medskip 

\noindent A. It is kind of spread out, which is why I wanted to write this note.

\begin{itemize}

\item These limit shapes are examples of Duistermaat-Heckman measures as defined by... well Duistermaat and Heckman \cite{DH}.  The general theory guarantees that these limits exist and are described by continuous, piecewise polynomial functions.

\item The specific example of Duistermaat-Heckman measures for $q$-binomials equaling a convolution of uniform measures on an interval is worked out as an example in \cite{CDKW}.

\item Stanley and Zanello studied the asymptotic behavior of the $\qbinom{n+k}{k}_q$ in \cite{SZ}.  They observed the local quasipolynomiality of the coefficients, and gave an explicit local approximation for the coefficients of $q^m$ with $m \approx \alpha nk$.  Suitably renormalized this can also be thought of as a local description of the limit shape. 

\end{itemize}

\end{document}